\documentclass[12pt,reqno,a4paper]{amsart}
\usepackage{amsmath,amssymb,amsthm,ifthen,color,amsaddr,a4wide}
\usepackage[all]{xypic}
\usepackage[mathcal]{euscript}
\usepackage{mathrsfs}

\newtheorem{thm}{Theorem}[section]

\newtheorem{prop}[thm]{Proposition}
\newtheorem{defin}[thm]{Definition}

\theoremstyle{remark}
\newtheorem{rem}[thm]{Remark}
\newtheorem{ex}[thm]{Example}

\newcommand{\cO}{{\mathcal O}}

\newcommand{\R}{{\mathbb R}}
\newcommand{\C}{{\mathbb C}}
\newcommand{\Z}{{\mathbb Z}}

\newcommand{\Q}{{\mathbb Q}}

\newcommand{\qee}{\mbox{\hspace{0.2mm}}\hfill$\triangle$}

\newcommand{\var}{{X_\Sigma}}

\newcommand{\marginnote}[1]{\ifthenelse{\isodd{\thepage}}{\normalmarginpar}
{\reversemarginpar}\marginpar{\fbox{\parbox{28mm}{\sloppy\footnotesize #1}}}}
\begin{document}
\begin{flushright} SISSA Preprint  20/2011/fm\end{flushright}
\bigskip
\title{$\Q$-factorial Laurent rings}
\author{\small Ugo Bruzzo$^{\S\dag\ddag}$ and Antonella Grassi$^\P$}
\address{\rm $^\S$ Institut des Hautes \'Etudes Scientifiques, Le Bois-Marie, \\
91440 Bures-sur-Yvette, France\footnote{$^\ddag$ On leave of absence from Scuola Internazionale Superiore di Studi Avanzati, Via Bonomea 265, 34136 Trieste, Italy}}
\address{\vskip-10pt\rm $^\dag$ Istituto Nazionale di Fisica Nucleare, Sezione di Trieste}
\address{\vskip-10pt\rm $^\P$ Department of Mathematics, University of Pennsylvania,\\
David Rittenhouse Laboratory, 209 S 33rd Street,\\ Philadelphia, PA 19104, USA\footnote{Support for this work was provided by the NSF Research Training Group Grant
DMS-0636606, by {\sc prin} ``Geometria delle variet\`a  algebriche
e dei loro spazi dei moduli '' and  the {\sc infn} project {\sc pi14} ``Nonperturbative dynamics of gauge theories''. U.B. is a member of the {\sc vbac} group.}}
 \thanks{E-mail: {\tt bruzzo@sissa.it, grassi@sas.upenn.edu}}
\date{\today}
\begin{abstract} Dolgachev
proves that   the ring naturally associated to a generic   Laurent polynomial in $d$ variables, $d \geq 4$, is  factorial  \cite{Dolg-Newton,Dolg-NewtonE} (for any field $k$).
We prove a sufficient condition for  the    ring   associated to a very general complex Laurent polynomial in $d=3$  variables to be $\Q$-factorial.
\end{abstract}

\maketitle

\section{Introduction}
In \cite{Dolg-Newton} and Dolgachev \cite{Dolg-NewtonE}  proves that the ring $A_F$ naturally associated to  generic   Laurent polynomial $F$ in $d$ variables, $d \geq 4$, with coefficients in any field $k$, is factorial. The basic ingredient in Dolgachev's proof
is Grothendieck's Lefschetz-type theorem (\cite{SGA2-XI}, Prop.~3.12) which, among other things, shows that under suitable conditions, the natural  restriction map $\operatorname{Pic}(X) \to \operatorname{Pic}(Y)$, where $X$ is a scheme and $Y$ is subvariety corresponding to an ideal sheaf in $\cO_X$, is an isomorphism. This result can be applied only when $d\ge 4$.

In this paper we consider the case $d=3$, assuming that $k=\C$, and prove a sufficient condition for the ring $A_F$ to be $\Q$-factorial (Theorem \ref{fact}). The proof of this fact follows the lines of Dolgachev's proof, with Grothendieck's result replaced by a Noether-Lefschetz theorem for hypersurfaces in toric 3-folds (Theorem \ref{BGthm}) that
we proved in \cite{BG}.

\bigskip

\noindent{\bf Acknowledgement.} We thank   Igor Dolgachev for interesting correspondence leading to this result and the referee for useful comments. The authors are grateful for the hospitality and support offered by the University of Pennsylvania, SISSA and IHES.

\bigskip

\section{Preliminaries}
We follow the notation in \cite{BaCox94} and \cite{BG}.
 Let $M$ be a $d$-dimensional lattice, $N=\operatorname{Hom}(M,\Z)$ and $\mathbf T_N= N \otimes \C^*$ the associated algebraic torus. Let $\Sigma\subset N_\R$ be a complete simplicial fan,
 and denote by $\var$ the corresponding complete toric variety. The torus $\mathbf  T_N$ naturally acts on $\var$;
$\mathbf T_\tau\subset \var$ denotes the orbit of a subset of $\var$ corresponding to a face
$\tau$ of $\Sigma$  under this action; the open dense orbit is denoted
by $\mathbf T_{0}$.

\begin{defin}\label{nondegenerate}\cite[Def.~4.13]{BaCox94}   A hypersurface $X$ in $\var$ is {\em nondegenerate} if $X\cap\mathbf T_\tau$ is a smooth 1-codimensional subvariety of
$\mathbf T_\tau$ for all faces $\tau$ in $\Sigma$.
\end{defin}

 $\var$ has only abelian quotient singularities, and is therefore an orbifold.

\begin{prop}\cite[Prop.~3.5, 4.15]{BaCox94}  Let $L$ be a ample line bundle on $\var$. The hypersurface $X\subset \var$ given by the zero locus of a generic section of $L$ is nondegenerate.
Moreover, $X$ is an orbifold.\label{nondeg}
\end{prop}

Since $X$ is an orbifold, its complex cohomology has a pure Hodge structure \cite{SaitoKyoto}.
This is an essential point in the proof of our Theorem \ref{BGthm}.

 \begin{defin}[The Cox Ring \cite{Cox95}] Consider a variable $z_i$ for each 1-dimensional cone $\varsigma_i$, $i=1,\dots,n$  in
 $\Sigma$, and let $S(\Sigma)$ be the polynomial ring $\C[z_1,\dots,z_n]$.
\end{defin}

The Cox ring has a natural gradation given by its class group $Cl(\Sigma)$ of $\var$.

Let $L$ be an ample line bundle on $\var$, and let $f\in H^0(\var,L)\simeq S(\Sigma)_\beta$,
where $\beta=\deg(L)$.

\begin{defin}  The
 {\em Jacobian ring}  of  $f$ is the quotient $R(f)=S(\Sigma)/J(f)$, where
$J(f)$ is  the ideal in $S(\Sigma)$ generated by the derivatives of $f$. \end{defin}

The Jacobian ring $R(f)$ inherits a natural gradation from $S(\Sigma)$.

The next theorem was proved in \cite{BG}, and will be key to proving our result about Laurent rings.
 We assume $d=3$. We recall that the Picard number is defined as the rank of the class group.

\begin{thm}\cite{BG} Let $\var$ a complete simplicial toric variety, and $X\subset \var$ a very general hypersurface cut by a section $f$
of an ample line bundle $L$  such that the multiplication morphism
\begin{equation}\label{Rf} R(f)_\beta\otimes R(f)_{\beta-\beta_0} \to R(f)_{2\beta-\beta_0} \end{equation}
is surjective (here $\beta=\deg(L)$ and $\beta_0=-\deg(K_\var)$, where $K_\var$ is the canonical sheaf of $\var$).
Then $X$ has the same Picard number as $\var$. 
\label{BGthm} \end{thm}
Recall that a property is very general if it holds in the
complement of countably many proper subvarieties.

If $X$ is a quartic  surface in $\mathbb P^3$, or more generally a $K3$ surface defined by a section of the anticanonical divisor in a simplicial  toric variety, then the above map is surjective \cite{BG}. It is a classical result that the map is not surjective if $X$ is a cubic in $\mathbb P^3$.

\bigskip

\section{$\Q$-factorial Laurent rings}

 The ring $\C[M]$ may be identified with the ring of regular functions on the torus $\mathbf T_N \simeq \mathbf T_{0}\subset \var$.
An element  $F\in \C[M]$  is called a {\em Laurent polynomial}; $F$ may be regarded as a section of the ample line bundle $L$, and it defines a hypersurface $X_F$ in $\var$. 

Let $\Delta\subset M\otimes_\Z\R$ be the polytope uniquely determined by the  fan $\Sigma$ and $L$ (see \cite{Oda88}, Lemma 2.14).  To each Laurent polynomial $F$ on can associate
a polytope $\Delta_F$, called the {\em Newton polytope} of $F$. This is most easily described by choosing an isomorphism $M\simeq \Z^d$, writing
$$ F = \sum_{i_1,\dots,i_d\in \Z^d} a_{i_1,\dots,i_d}\,t_1^{i_1}\cdots t_d^{i_d}$$
and defining
$$\operatorname{supp}(F)= \{i_1,\dots,i_d\in \Z^d\,\vert a_{i_1,\dots,i_d}\ne 0 \}.$$
 $\Delta_F$ is then defined to be the  convex hull   of $\operatorname{supp}(F)$ and $\Gamma(\Delta)$  the set of all Laurent polynomials such that $\Delta_F\subset\Delta$. $\Gamma(\Delta)$  is a finite dimensional vector space over $\C$.

 By results given in \cite{Kov78} (see also \cite{Oda88}, Chapter 2) a Laurent polynomial $F$
 extends to a meromorphic function on $\var$, which is
 a section  of an ample line bundle $L_F$. Thus, $F$ may be regarded as an element in $S(\Sigma)_\beta$,
where $\beta=\deg(L_F)$. Denote by $A_F$   the ring $\C[M]/(F)$.

\begin{thm} \label{fact} Let $d=3$,  and let $F$ be a very general   Laurent polynomial in $\Gamma(\Delta)$; set $\beta=\deg(L_F)$ and $\beta_0=-\deg(K_\var)$.
 If  the multiplication morphism
\begin{equation}\label{RF} R(F)_\beta\otimes R(F)_{\beta-\beta_0} \to R(F)_{2\beta-\beta_0} \end{equation}
 is surjective, the ring $A_F$  is $\Q$-factorial.
\end{thm}

The proof that $A_F$ is $\Q$-factorial follows closely the proof of Theorem 1.1 in \cite{Dolg-Newton}. The basic idea is to formulate   the problem in a geometric way:
\begin{proof}Let  $X_F \subset \var$ be the hypersurface  cut by $F$ (as a section of $L_F$).
By Proposition \ref{nondeg} the hypersurface $X_F$ is nondegenerate, and is an orbifold.

 Note that the ring $A_F$ may be identified with the ring of regular functions on the affine part $U_F=X_F\cap \mathbf T_{0}$ of $X_F$. Since the Picard group of $\mathbf T_{0}$ is trivial,
 every Cartier  divisor in $\var$ is linearly equivalent to a divisor supported in  $\var-\mathbf T_{0}$.
By Theorem \ref{BGthm}, $X_F$ has the same Picard number as $\var$, i.e., $\rho(X_F)=\rho(\var)$.
 Then any Cartier divisor in $X_F$ is linearly equivalent modulo torsion to a divisor supported in $X_F-U_F$,  so that $\operatorname{Pic}(U_F)\otimes\Q =0$. Since $U_F$ is normal (actually smooth), then
 $Cl(U_F)\otimes\Q =0$. As
$U_F\simeq \operatorname{Spec}(A_F)$, we have $Cl(A_F)\otimes \Q=0$.
\end{proof}

\bigskip

 \bigskip\frenchspacing

\def\cprime{$'$}

\end{document}